\documentclass[12pt]{article}

\usepackage{amscd,amsmath, amssymb, fancyhdr, epsfig,color,enumitem}
\definecolor{dblue}{rgb}{0,0,.6}
\usepackage[colorlinks=true, linkcolor=dblue, citecolor=dblue, filecolor = dblue, menucolor = dblue, urlcolor = dblue]{hyperref}

\numberwithin{equation}{section}

\newcommand{\version}{version 3.2,\ \   Febr. 8, 2019}

\def\eqref#1{(\ref{#1})}

\newcommand{\arrow}{{\:\longrightarrow\:}}
\newcommand{\Z}{{\Bbb Z}}
\newcommand{\C}{{\Bbb C}}

\newcommand{\R}{{\Bbb R}}

\def\1{\sqrt{-1}\:}

\newcommand{\cntrct}                
{\hspace{2pt}\raisebox{1pt}{\text{$\lrcorner$}}\hspace{2pt}}

\makeatletter
\def\x@arrow{\DOTSB\Relbar}
\def\xlongequalsignfill@{\arrowfill@\x@arrow\Relbar\x@arrow}
\newcommand{\xlongequal}[2][]{%
        \ext@arrow 0099\xlongequalsignfill@{#1}{#2}}
\def\xlongrightarrowfill@{\arrowfill@\relbar\relbar\longrightarrow}
\newcommand{\xlongrightarrow}[2][]{%
        \ext@arrow 0099\xlongrightarrowfill@{#1}{#2}}
\makeatother

\renewcommand{\bar}{\overline}
\renewcommand{\phi}{\varphi}
\renewcommand{\epsilon}{\varepsilon}
\renewcommand{\geq}{\geqslant}
\renewcommand{\leq}{\leqslant}

\newcommand{\Sp}{\operatorname{Sp}}
\newcommand{\End}{\operatorname{End}}
\newcommand{\Mat}{\operatorname{Mat}}

\newcommand{\Cl}{\operatorname{{\mathcal Cl}}}

\newcommand{\Hol}{\operatorname{Hol}}
\newcommand{\Aut}{\operatorname{Aut}}

\newcommand{\Sec}{\operatorname{Sec}}

\newcommand{\Diff}{\operatorname{\sf Diff}}

\newcommand{\Tw}{\operatorname{Tw}}

\newcommand{\Spin}{\operatorname{Spin}}

\newcommand{\Teich}{\operatorname{\sf Teich}}
\newcommand{\Comp}{\operatorname{\sf Comp}}

\newcommand{\bbZ}{\mathbb{Z}}

\newcommand{\bbR}{\mathbb{R}}
\newcommand{\bbC}{\mathbb{C}}

\newcommand{\bbP}{\mathbb{P}}

\newcommand{\emrp}{\mathrm{End}}

\newcommand{\Pin}{\mathrm{Pin}}

\newcommand{\bOmega}{{\boldsymbol{\boldsymbol{\Omega}}}}

\newcounter{Mycounter}[section]
\newcounter{lemma}[section]
\renewcommand{\thelemma}{{Lemma \thesection.\arabic{lemma}}}
\newcommand{\lemma}{%
    \setcounter{lemma}{\value{Mycounter}}
    \refstepcounter{lemma}
    \stepcounter{Mycounter}
    {\noindent \bf \thelemma:\ }}

\newcounter{claim}[section]
\renewcommand{\theclaim}{{Claim \thesection.\arabic{claim}}}
\newcommand{\claim}{%
    \setcounter{claim}{\value{Mycounter}}
    \refstepcounter{claim}
    \stepcounter{Mycounter}
    {\noindent \bf \theclaim:\ }}

\newcounter{sublemma}[section]
\newcounter{corollary}[section]
\renewcommand{\thecorollary}{{Corollary \thesection.\arabic{corollary}}}
\newcommand{\corollary}{%
    \setcounter{corollary}{\value{Mycounter}}
    \refstepcounter{corollary}
    \stepcounter{Mycounter}
    {\noindent \bf \thecorollary:\ }}

\newcounter{theorem}[section]
\renewcommand{\thetheorem}{{Theorem \thesection.\arabic{theorem}}}
\newcommand{\theorem}{%
    \setcounter{theorem}{\value{Mycounter}}
    \refstepcounter{theorem}
    \stepcounter{Mycounter}
    {\noindent \bf \thetheorem:\ }}

\newcounter{conjecture}[section]
\newcounter{proposition}[section]
\renewcommand{\theproposition}
      {{Proposition \thesection.\arabic{proposition}}}
\newcommand{\proposition}{%
    \setcounter{proposition}{\value{Mycounter}}
    \refstepcounter{proposition}
    \stepcounter{Mycounter}
    {\noindent \bf \theproposition:\ }}

\newcounter{definition}[section]
\renewcommand{\thedefinition}
      {{Definition~\thesection.\arabic{definition}}}
\newcommand{\definition}{%
    \setcounter{definition}{\value{Mycounter}}
    \refstepcounter{definition}
    \stepcounter{Mycounter}
    {\noindent \bf \thedefinition:\ }}

\newcounter{example}[section]
\renewcommand{\theexample}{{Example \thesection.\arabic{example}}}
\newcommand{\example}{%
    \setcounter{example}{\value{Mycounter}}
    \refstepcounter{example}
    \stepcounter{Mycounter}
    {\noindent \bf \theexample:\ }}

\newcounter{remark}[section]
\renewcommand{\theremark}{{Remark \thesection.\arabic{remark}}}
\newcommand{\remark}{%
    \setcounter{remark}{\value{Mycounter}}
    \refstepcounter{remark}
    \stepcounter{Mycounter}
    {\noindent \bf \theremark:\ }}

\newcounter{problem}[section]
\newcounter{question}[section]
\renewcommand{\thequestion}{{Question \thesection.\arabic{question}}}
\newcommand{\question}{%
    \setcounter{question}{\value{Mycounter}}
    \refstepcounter{question}
    \stepcounter{Mycounter}
    {\noindent \bf \thequestion:\ }}

\newcommand{\proof}{{\bf Proof:\:}}

\makeatletter

\renewcommand{\leftmark}{{\scriptsize  
$k$-symplectic structures and absolutely trianalytic subvarieties}}

\@addtoreset{equation}{section} \@addtoreset{footnote}{section} \makeatother

\def\blacksquare{\hbox{\vrule width 5pt height 5pt depth 0pt}}
\def\endproof{\blacksquare}

\begin{document}
\begin{center}
{\LARGE\bf
$k$-symplectic structures and absolutely trianalytic
subvarieties in hyperk\"ahler manifolds.\\[3mm]
}

Andrey Soldatenkov\footnote{Andrey Soldatenkov is partially supported by
AG Laboratory NRU-HSE, RF government grant, ag. 11.G34.31.0023,
MK-1297.2014.1 and a grant from Dmitri Zimin's ``Dynasty'' foundation},
Misha Verbitsky\footnote{Partially supported by RSCF grant 14-21-00053
within AG Laboratory NRU-HSE.}

\end{center}

{\small \hspace{0.02\linewidth}
\begin{minipage}[t]{0.85\linewidth}
\centerline{{\bf Abstract}}
Let $(M,I,J,K)$ be a hyperk\"ahler manifold, and $Z\subset (M,I)$
a complex subvariety in $(M,I)$. We say that $Z$ is trianalytic if
it is complex analytic with respect to $J$ and $K$, and absolutely
trianalytic if it is trianalytic with respect to any hyperk\"ahler
triple of complex structures $(M,I,J',K')$ containing $I$.
For a generic complex structure $I$ on $M$, all complex 
subvarieties of $(M,I)$ are absolutely trianalytic.
It is known that the normalization $Z'$ of a trianalytic
subvariety is smooth; we prove that $b_2(Z')\geq b_2(M)$,
when $M$ has maximal holonomy (that is, $M$ is IHS).

\qquad To study absolutely trianalytic subvarieties further, we 
define a new geometric structure, called $k$-symplectic 
structure; this structure is a generalization of 
hypersymplectic structure. A $k$-symplectic structure
on a $2d$-dimensional manifold $X$ is a $k$-dimensional space $R$ 
of closed 2-forms on $X$ which all have rank $2d$ or $d$. 
It is called non-degenerate if the set of all
degenerate forms in $R$ is a smooth, non-degenerate
quadric hypersurface in $R$.  

\qquad We consider absolutely trianalytic tori
in a hyperk\"ahler manifold $M$ of maximal holonomy. 
We prove that any such torus is equipped with a 
non-degenerate $k$-symplectic 
structure, where $k=b_2(M)$. We show that the tangent
bundle $TX$ of a $k$-symplectic manifold is a Clifford
module over a Clifford algebra $Cl(k-1)$.
Then an absolutely trianalytic
torus in a hyperk\"ahler manifold $M$ with 
$b_2(M)\geq 2r+1$ is at least $2^{r-1}$-dimensional.
\end{minipage}
}

\tableofcontents


\section{Introduction}


\subsection{Absolutely trianalytic subvarieties in
  hyperk\"ahler manifolds}

Let $M$ be a K\"ahler, compact, holomorphic symplectic
manifold. Calabi-Yau theorem (\cite{_Yau:Calabi-Yau_})
implies that $M$ admits a Ricci-flat metric $g$, unique
in each K\"ahler class. Using Berger's classification of
Riemannian holonomies and Bochner vanishing, one shows
that the Levi-Civita connection of $g$ preserves a triple
of complex structures $I,J,K$ satisfying the quaternionic
relation $IJ =-JI=K$ (\cite{_Besse:Einst_Manifo_}).
A Riemannian manifold admitting a triple of complex structures
$I,J,K$ satisfying quaternionic relations and K\"ahler with
respect to $g$ is called {\bf hyperk\"ahler}. One can construct a
holomorphic symplectic form on any hyperk\"ahler manifold as follows.
There are K\"ahler forms $\omega_I$, $\omega_J$ and $\omega_K$ associated
to the complex structures $I,J$ and $K$. One can check that the
2-form $\omega_J+\sqrt{-1}\omega_K$ is of Hodge type $(2,0)$ with
respect to $I$. Since it is also closed it is a holomorphic symplectic form.
So we shall treat the terms ``hyperk\"ahler'' and ``holomorphic
symplectic'' as (essentially) synonyms.

Given any triple $a, b ,c \in \R$, $a^2+b^2+c^2=1$,
the operator $L=aI+bJ+cK$ satisfies $L^2=-1$ and defines 
a K\"ahler structure on $(M,g)$; such a complex structure is called
{\bf induced by the hyperk\"ahler structure}. Complex subvarieties
of such $(M,L)$ for generic $(a,b,c)$ were studied in 
\cite{_Verbitsky:trianaly1_}, 
\cite{_Verbitsky:trianaly2_} and \cite{_Verbitsky:subvar_}.

\hfill

\definition
Let $Z\subset M$ be a closed subset of a 
hyperk\"ahler manifold. It is called {\bf trianalytic},
if it is complex analytic with respect to all induced
complex structures $L$.

\hfill

In the definition above it is enough to require the
subvariety $Z$ to be complex analytic with respect only to $I$ and $J$.
Then it will automatically be complex analytic with respect to any
induced complex structure.
This is clear, because $Z$ is trianalytic if and only if
for all smooth points $z\in Z$, the space
$T_zZ\subset T_zM$ is preserved by the quaternion algebra ${\Bbb H}$
(\cite{_Verbitsky:trianaly2_}).
However, ${\Bbb H}$ is generated by any two non-collinear 
elements $I,I_1$ with $I^2=I_1^2=-1$.

\hfill

Singularities of trianalytic subvarieties always admit
a hyperk\"ahler resolution.

\hfill

\theorem\label{_Verbitsky_desing_Theorem_} (\cite{_Verbitsky:desing_})
Let $M$ be a hyperk\"ahler manifold, $Z\subset M$ a trianalytic
subvariety, and $I$ an induced complex structure.
Consider the normalization \[ \widetilde{(Z, I)}\stackrel n\arrow (Z,I)\] 
of $(Z,I)$. Then $\widetilde{(Z, I)}$ is smooth, and
the map $\widetilde{(Z, I)}\arrow M$ is an immersion,
inducing a hyperk\"ahler structure on $\widetilde{(Z, I)}$.
\endproof

\hfill

\definition\label{_twistor_family_Definition_}
Let $M$ be a hyperk\"ahler manifold, and $S$ the family
of all induced complex structures $L=aI+bJ+cK$, where
$a, b ,c \in \R$, $a^2+b^2+c^2=1$. Then $S$ is called
{\bf the twistor family} of complex structures.

\hfill

The following theorem implies that whenever $L$ is a generic
element of a twistor family, all subvarieties
of $(M,L)$ are trianalytic.

\hfill

\theorem (\cite{_Verbitsky:trianaly2_,_Verbitsky:subvar_}) 
Let $M$ be a hyperk\"ahler manifold, $S$ its twistor family.
Then there exists a countable subset $S_1\subset S$, such that
for any complex structure $L\in S \backslash S_1$, all
compact complex subvarieties of $(M,L)$ are trianalytic.
\endproof

\hfill

In \cite{_Verbitsky:Hilbert_} 
(see also \cite{_KV:book_V_}, \cite{_KV:part-res_} 
and \cite{_Oguiso:automo_}),
this theorem was used to study subvarieties of generic
deformations of a compact holomorphic symplectic manifold $M$.
Recall that the Teichm\"uller space of $M$ is the quotient
$\Teich = \Comp/\Diff^0$ of the (infinite-dimensional) space of
all complex structures of hyperk\"ahler type by the group
$\Diff^0$ of isotopies (\cite{_V:Torelli_}).
$\Teich$ is a complex, non-Hausdorff manifold. If we fix a complex structure
$I$ on $M$ then the connected component of the Teichm\"uller space containing
$I$ can be identified with a connected component of the so-called ``marked moduli space''
of deformations of $(M, I)$ (\cite{_V:Torelli_}).

\hfill

\definition\label{_abs_triana_Definition_}
Let $(M,I,J,K)$ be a compact, holomorphic symplectic,
K\"ahler manifold, and $Z\subset (M,I)$ a complex subvariety,
which is trianalytic with respect to any hyperk\"ahler structure
compatible with $I$. Then $Z$ is called {\bf absolutely trianalytic}.

\hfill

\definition
For a given complex structure $I$, consider {\bf the Weil operator}
$W_I$ acting on $(p,q)$ forms as $\1(p-q)$.
Let $G_{MT}(M,I)$ be a smallest rational algebraic 
subgroup of $\Aut(H^*(M,\R))$ containing $e^{tW_I}$.
This group is called {\bf the Mumford-Tate group of $(M,I)$}.
A group generated by $G_{MT}(M,I)$ for all complex structures 
$I$ in a connected component of a deformation space is called 
{\bf a maximal Mumford-Tate group of $M$} (\cite{_Deligne:Hodge_cycles_}).

\hfill

It is not hard to check that the 
Mumford-Tate group $G_{MT}(I)$ of $(M,I)$ is lower semicontinuous
as a function of $I\subset \Teich$ in Zariski topology
on $\Teich$ (\cite{_Deligne:Hodge_cycles_}).
This implies that $G_{MT}(I)$ is constant outside
of countably many complex subvarieties of positive
codimension. We call $I\in \Teich$ {\bf Mumford-Tate generic}
if $G_{MT}(I)$ is maximal. If $M$ has maximal holonomy, 
the maximal Mumford-Tate group is isomorphic to
$\Spin(H^2(M,\R),q)$ (\cite{_KV:book_V_}).
Any $I\in \Teich$ outside of a countably many subvarieties
of positive codimension is Mumford-Tate generic.

\hfill

\remark\label{_MT_gene_defo_Remark_}
Let $I$ be a Mumford-Tate generic complex structure,
and $\eta$ an integer $(p,p)$-class. Then $\eta$
is of type $(p,p)$ for any deformation of $I$.

\hfill

Absolutely trianalytic subvarieties can be characterized
in terms of the Mumford-Tate group, as follows.

\hfill

\claim\label{_gene_MT_abs_triana_Claim_}
Let $(M,I,J,K)$ be a hyperk\"ahler manifold, and $Z\subset (M,I)$
be a complex subvariety. Then $Z$ is absolutely trianalytic
if and only if its fundamental class is $G$-invariant,
where $G$ is a maximal Mumford-Tate group of $M$.
In particular, $Z$ is absolutely trianalytic
when $(M,I)$ is Mumford-Tate generic.

\hfill

{\bf Proof:} This statement  
follows from the definitions (\cite[Claim 4.4]{_V:cohe_K3_}).
\endproof

\hfill

Clearly, the set of absolutely trianalytic subvarieties
does not change if one passes from one complex structure
in a twistor family $S$ to another complex structure in $S$.
On the other hand, any two complex structures in the same
component of $\Teich$ can be connected by a sequence of twistor
families (\cite{_Verbitsky:coho_announce_}). 
This means that the set of absolutely trianalytic subvarieties
in $(M, I)$ is determined by the connected component of
$\Teich$ where $I$ lies.

This can be used to prove the following theorem.

\hfill

\theorem \label{_abs_triana_diffeo_Theorem_}
Let $I_1, I_2\in \Teich$ be points in the same connected component
of the Teichm\"uller space. Then there exists a diffeomorphism
$\nu:\; (M,I_1)\arrow (M, I_2)$ such that any absolutely trianalytic
subvariety $Z\subset (M,I_1)$ is mapped to an absolutely trianalytic
subvariety $\nu(Z)\subset (M, I_2)$.
\endproof

\hfill

Absolutely trianalytic subvarieties were studied 
in \cite{_Verbitsky:Hilbert_}, where it was shown that
a general deformation of a Hilbert scheme of a K3 surface
has no complex (or, equivalently, no absolutely trianalytic)
subvarieties.  In \cite{_KV:gen_Kum_}, Kaledin and Verbitsky
used the same argument to study absolutely tri\-analytic
subvarieties in generalized Kummer varieties. They 
``proved'' non-existence of such subvarieties, but their 
argument was faulty (in fact, there exists
an absolutely trianalytic subvariety in 
a generalized Kummer variety). In \cite{_KV:part-res_}, 
the error was found, and the argument was repaired
to show that any absolutely trianalytic subvariety 
$Z$ of a generalized Kummer
variety is a deformation of a resolution of 
singularities of a quotient of a torus by a Weyl group action.
The idea was to show that $Z$ is a resolution of 
singularities of a quotient of a flat subtorus in a symmetric
power of a torus, and classify the group actions which 
admit a holomorphic symplectic resolution.
Later, Ginzburg and Kaledin have shown that only the Weyl groups $A_n$,
$B_n$, $C_n$ can occur in quotient maps with
quotients which admit holomorphically symplectic
resolution of singularities (\cite{_GinzKale:Poisson_}).

Non-existence of absolutely trianalytic subvarieties
in a Hilbert scheme $M$ of K3 was used in 
\cite{_KV:book_V_} to prove compactness of 
deformation spaces of certain stable holomorphic bundles
on $M$.

\subsection{Hyperk\"ahler manifolds}

Before we state our main results, let us introduce the
Bogomolov-Beauville-Fujiki form and maximal holonomy manifolds.

\hfill

\definition \label{_IHS_Definition_}
A compact hyperk\"ahler manifold $M$ is called
{\bf simple}, or {\bf maximal holonomy}, or {\bf IHS}
(from ``Irreducible Holomorphic Symplectic'')  
if $\pi_1(M)=0$ and $H^{2,0}(M)=\C$.

\hfill

\remark 
It follows from Bochner's vanishing
and Berger's classification of holonomy groups
that a hyperk\"ahler manifold has maximal holonomy 
$\Sp(n)$ whenever $\pi_1(M)=0$, $H^{2,0}(M)=\C$
(\cite{_Besse:Einst_Manifo_}).
This explains the term.

\hfill

\theorem \label{_Bogo_decompo_Theorem_}
(Bogomolov's decomposition; \cite{_Bogomolov:decompo_})
Any compact hyperk\"ahler manifold admits a finite covering
which is a product of a torus and several simple hyperk\"ahler manifolds.
\endproof

\hfill

\theorem \label{_Fujiki_formula_Theorem_}
(Fujiki, \cite{_Fujiki:HK_})
Let $M$ be a simple hyperk\"ahler
manifold of complex dimension $2n$. Then there exists
a primitive integral quadratic form $q$ on $H^2(M, \Z)$ and
a constant $c_M$ such that for any $\eta\in H^2(M, \C)$ we have
$\int_M \eta^{2n}=c_M q(\eta,\eta)^n$.
\endproof

\hfill

\remark\label{_BBF_Remark_}
\ref{_Fujiki_formula_Theorem_} determines the form $q$
uniquely up to a sign. For $\dim_\R M=4n$,
$n$ odd, sign is also determined. For $n$ even, the sign is determined by
the following formula, due to Bogomolov and Beauville.
\begin{equation}\label{_BBF_expli_Equation_}
\begin{aligned} 
\mu q(\eta,\eta) &=
   (n/2)\int_X \eta\wedge\eta  \wedge \Omega^{n-1}
   \wedge \bar \Omega^{n-1} -\\
 &-(1-n)\left(\int_X \eta \wedge \Omega^{n-1}\wedge \bar
   \Omega^{n}\right) \left(\int_X \eta \wedge \Omega^{n}\wedge \bar \Omega^{n-1}\right)
\end{aligned}
\end{equation}
where $\Omega$ is the holomorphic symplectic form, and 
$\mu>0$ a positive constant.

\hfill

\definition
Let $M$ be a hyperk\"ahler manifold of maximal holonomy,
and $q$ the form on $H^2(M)$ defined by \ref{_BBF_Remark_} and 
\ref{_Fujiki_formula_Theorem_}. Then $q$ is called 
{\bf Bogomolov-Beauville-Fujiki form} (BBF form).
The equation \eqref{_BBF_expli_Equation_} (or, even better,
a similar equation \cite[(1.1)]{_Verbitsky:cohomo_},
 using the K\"ahler form instead
of the holomorphic symplectic form) can be used
to show that $q$ is non-degenerate and has signature $(3, b_2(M)-3)$.

\hfill

We will need Fujiki relations in
greater generality which we also recall 
(the proof can be found in \cite{GHJ}). Namely, let $M$ be a simple hyperk\"ahler manifold
with $\dim_\C M = 2n$. Consider
a cohomology class $\gamma$ which is of Hodge type $(2n-2m,2n-2m)$ on any small deformation
of $M$.
Then there exists a constant $c_\gamma\in \bbR$ such that for any $\alpha, \beta\in H^2(X, \bbC)$ we have
\begin{equation}\label{_pairing_subvari_Equation_}
\gamma\cdot \alpha^{2m-1}\cdot \beta = c_\gamma q(\alpha,\alpha)^{m-1}q(\alpha,\beta).
\end{equation}
In this equation on the left hand side we have the intersection product in cohomology of $M$.

\subsection{Main results: $k$-symplectic structures 
and absolutely trianalytic subvarieties}

In the present paper, we prove two bounds on the Betti 
numbers of absolutely trianalytic subvarieties, quite restrictive
for their geometry. In particular, we prove that 
the normalization of a proper complex subvariety
of a generic deformation of a 10-dimensional
O'Grady space must necessarily belong  to a new type
of hyperk\"ahler manifolds (\ref{_O_Gra_10_dim_Theorem_}).

Our arguments are based on an elementary observation, stated below
as \ref{_form_restri_Theorem_}.

Recall (\ref{_Verbitsky_desing_Theorem_}) that the normalization of any trianalytic subvariety $Z\hookrightarrow M$
is a smooth hyperk\"ahler manifold $\tilde Z$ immersed into $M$.
Therefore, we can replace any trianalytic cycle by an immersed
hyperk\"ahler manifold. Note that the complex dimension of any
trianalytic subvariety is even.

\hfill

\theorem\label{_form_restri_Theorem_}
Let $M$ be a maximal holonomy hyperk\"ahler manifold and
$\tilde Z\stackrel\phi \arrow M$ 
the normalization of an absolutely trianalytic cycle $Z\subset M$
or a finite covering of such normalization.
Consider the polynomial $P_{\tilde Z}$ on $H^2(\tilde Z, \C)$
mapping a cohomology class $\eta$ to
$\int_{\tilde Z} \eta^{\dim_\C {\tilde Z}}$. Then for any $\alpha\in H^2(M,\C)$,
one has $P_{\tilde Z}(\phi^*\alpha)=\deg(\phi) c_Zq(\alpha,\alpha)^{\frac 1 2 \dim_\C Z}$,
where $c_Z$ is a constant determined by $M$ and $Z$.

\hfill

\proof Denote by $\gamma$ the fundamental class of $Z$ in the cohomology
of $M$. We have $P_{\tilde Z}(\phi^*\alpha) = \deg(\phi)\gamma\cdot\alpha^{\dim_\C Z} = c_Zq(\alpha,\alpha)^{\frac 1 2 \dim_\C Z}$.
The last equality follows from \eqref{_pairing_subvari_Equation_} where we put $\beta=\alpha$ and $c_Z=c_\gamma$. \endproof

\hfill

This simple result has very nice consequences.

\hfill

\corollary\label{_inje_h^2_Corollary_}
Let $M$ be a maximal holonomy hyperk\"ahler manifold,
and $\tilde Z\stackrel\phi \arrow M$ the normalization
of an absolutely trianalytic cycle $Z\subset M$,
or a finite covering of the normalization.
Then the induced map $H^2(M,\C) \stackrel {\phi^*}\arrow H^2(\tilde Z, \C)$
is injective.

\hfill

{\bf Proof:} Given $x\in \ker \phi^*\subset H^2(M,\C)$
and any $y\in H^2(M, \C)$, one would obtain
\begin{eqnarray}
\deg(\phi)c_Zq(x+y,x+y)^{\frac 1 2 \dim_\C Z}&=&P_{\tilde Z}(\phi^*(x+y))= P_{\tilde Z}(\phi^*y)=\nonumber\\
 &=& \deg(\phi)c_Zq(y,y)^{\frac 1 2 \dim_\C Z}.\nonumber
\end{eqnarray}
This gives $q(x,y)=0$ for all $y\in H^2(M, \C)$.
However, this implies $x=0$, because $q$ is non-degenerate.
\endproof

\hfill

%


\corollary\label{_b_2_bound_abs_triana_IHS_Corollary_}
Let $M$ be a maximal holonomy hyperk\"ahler manifold,
$Z\subset M$ an absolutely trianalytic subvariety, and
$\tilde Z\arrow M$ its normalization. Consider a 
finite covering $\tilde Z_1\arrow \tilde Z$
such that $\tilde Z_1$ is a product of a hyperk\"ahler
torus $T$ and several maximal holonomy hyperk\"ahler manifolds $K_i$
(such a decomposition always exists by
Bogomolov decomposition \ref{_Bogo_decompo_Theorem_}).
Then $b_2(T)\geq b_2(M)$ and $b_2(K_i)\geq b_2(M)$.

\hfill

{\bf Proof:} Any hyperk\"ahler metric on $M$ induces a hyperk\"ahler
structure on $\tilde Z_1$. Each component in Bogomolov decomposition
of $\tilde Z_1$ is immersed into $M$ as a hyperk\"ahler subvariety.
Hence the images of those components are also absolutely trianalytic subvarieties
of $M$ and we can apply \ref{_inje_h^2_Corollary_} to them.
\endproof

\hfill

In this paper we are also interested in absolutely trianalytic
subvarieties $Z\subset M$ such that the normalization $\tilde Z$
is a torus. We prove the following theorem.

\hfill

\theorem\label{_tori_main_bound_Theorem_}
Let $M$ be a hyperk\"ahler manifold of maximal holonomy, 
$T$ a hyperk\"ahler torus, and $T\arrow M$ a hyperk\"ahler
immersion with absolutely trianalytic image. Then 
$b_1(T) \geqslant 2^{\left\lfloor(b_2(X)-1)/2\right\rfloor}.$

{\bf Proof:} \ref{_b_1_torus_Corollary_}. \endproof

\hfill

Together with \ref{_b_2_bound_abs_triana_IHS_Corollary_},
this allows to prove non-existence of subvarieties
of known type in a 10-dimensional O'Grady manifold $M$
(\ref{_10_dim_subva_Remark_}), which has $b_2(M)=24$. 

\hfill

\theorem\label{_O_Gra_10_dim_Theorem_}
Let $Z\subset M$ be a proper complex subvariety of a general
deformation of 10-dimensional O'Grady manifold,
$\tilde Z$ its normalization, and $\tilde Z_1$
its covering equipped with the Bogomolov decomposition obtained as in
\ref{_b_2_bound_abs_triana_IHS_Corollary_}.
Then $\tilde Z_i = \prod_i K_i$ where $K_i$ are
maximal holonomy hyperk\"ahler manifolds
with $b_2 \geq 24$, that is, previously unknown type.
\endproof

\hfill

\remark
A maximal holonomy hyperk\"ahler manifold 
$K_i$ with $\dim_\C K_i=2$ satisfies $b_2(K_i) \leq 22$  
(Kodaira-Enriques classification, see
e.g. \cite{_Besse:Einst_Manifo_}). When $\dim_\C K_i=4$, one has
$b_2(K_i) \leq 23$ (\cite{_Guan:dim8_}). Therefore,
any absolutely trianalytic subvariety in a 10-dimensional O'Grady manifold
(if it exists) satisfies $\dim_\C Z\geq 6$, and has maximal holonomy.

\hfill

For a 6-dimensional O'Grady manifold $M$, one has
$b_2(M)=8$, so we can not prove analogous result by our methods.
However, complex tori are still prohibited.

\hfill

\theorem\label{_gene_defo_6_dim_Ogda_no_tori_Theorem_}
Let $M$ be a  general
deformation of a 6-dimensional O'Grady manifold. Then any
holomorphic map from a complex torus to $M$ is trivial.

\hfill

{\bf Proof:} 
Since an image of such a map is absolutely trianalytic,
its normalization is hyperk\"ahler. However, a dominant holomorphic
map from a torus $T$ to a manifold $Z$ with trivial canonical
bundle is a projection to a quotient torus. This is clear
because fibers of such map have trivial canonical bundle
by adjunction formula, but a Calabi-Yau submanifold
in a torus is also a torus. Then \ref{_gene_defo_6_dim_Ogda_no_tori_Theorem_}
 follows from 
\ref{_6_dim_subva_Remark_}. \endproof


\section{$k$-symplectic structures and Clifford representations}


The proof of \ref{_tori_main_bound_Theorem_}
is based on a discovery of a previously unknown geometric
structure, called $k$-symplectic structure. In 
Subsection \ref{_k_symple_on_vector_Subsection_} we study
$k$-symplectic structures on vector spaces, and in
Subsection \ref{_k_symple_appli_Subsection_}
we apply these linear-algebraic results to 
algebraic geometry of absolutely trianalytic tori.

$k$-symplectic structures generalize the hypersymplectic
structures known for a long time 
(\cite{_Arnold:hypersy_}, \cite{_Dancer_Swann:hypersy_}, 
\cite{_Andrada_Dotti_}), and trisymplectic
(3-symplectic) structures  defined in 
\cite{_Jardim_Verbitsky:trisymple_}.

\subsection{$k$-symplectic structures on vector spaces}
\label{_k_symple_on_vector_Subsection_}

Consider a complex vector space $V$ of dimension
$\dim_\bbC V = 4n$. Let $k$ be a non-negative integer.

\hfill

\definition
  A $k$-symplectic structure on $V$ is a subspace $\bOmega\subset \Lambda^2 V^*$ of dimension $k$,
  such that for some non-zero quadratic form $q\in S^2 \bOmega^*$ the following condition is satisfied:
  for any non-zero $\omega\in \bOmega$ we have
  $$
  \dim(\ker\omega) = 
  \begin{cases}
    2n,&\text{if $q(\omega) = 0;$}\\
    0,&\text{otherwise.}
  \end{cases}
  $$
  A $k$-symplectic structure is called non-degenerate if 
the quadratic form $q$ is non-degenerate.
  A vector space with a $k$-symplectic structure will be called a $k$-symplectic vector space.

\hfill

  A few remarks about this definition.

\hfill

\remark
The quadratic form $q$ in the definition above is unique up to a non-zero multiplier.
    We can reformulate the condition from the definition as follows. Consider
    the Pfaffian variety $$\mathcal P = \{\omega\in \bbP(\Lambda^2V^*)\mid \omega^{2n}=0\}$$ in $\bbP(\Lambda^2V^*)$.
    A subspace $\bOmega\subset\Lambda^2V^*$ defines a $k$-symplectic structure if and only
    if $\bbP\bOmega\cap\mathcal P$ is a quadric (necessarily of multiplicity $n$),
    and all the forms $\omega$ lying on this quadric have rank $2n$.
    The quadratic form $q$ is the one defining this quadric. If it is necessary to mention
    the form $q$ explicitly, we will denote a $k$-symplectic structure by $(\bOmega, q)$.

\hfill

\remark
 A $1$-symplectic structure is the same thing as a non-zero two-form $\omega\in \Lambda^2 V^*$
    which is either non-degenerate or has rank $2n$.

\hfill

\remark
Consider a non-degenerate $2$-symplectic structure $(\bOmega, q)$ on $V$.
    The quadric in $\bbP\bOmega=\bbP^1$ defined by $q$ is non-degenerate, hence it
    consists of two distinct points. Denote the corresponding two-forms by $\omega_1$ and $\omega_2$.
    The kernels $V_1$ and $V_2$ of these forms have dimension $2n$ by definition of a 2-symplectic structure.
    Moreover, $V_1\cap V_2 = 0$, because the generic linear combination of $\omega_1$ and $\omega_2$
    must be non-degenerate. Then we have $V = V_1\oplus V_2$.

    Conversely, given a symplectic vector space $(W,\omega)$, consider the direct sum $V = W^{\oplus2}$
    and denote by $\pi_i\colon V\to W$ the projection to the $i$-th summand. Then the subspace
    $\bOmega\subset \Lambda^2V^*$ spanned by $\pi_1^*\omega$ and $\pi_2^*\omega$ defines a
    non-degenerate 2-symplectic structure on $V$.

\hfill



\remark
Given a $k$-symplectic structure $(\bOmega, q)$, we can consider any subspace
    $\bOmega'\subset \bOmega$. Then
    $(\bOmega',q|_{\bOmega'})$ is a $k'$-symplectic structure for $k'=\dim_\bbC \bOmega'$.
    We will call it a substructure
    of $(\bOmega,q)$. Note that such substructure can be degenerate even if the initial structure was not.

\hfill

\remark
Consider a vector space $V$ of complex dimension 4. Then $\bOmega = \Lambda^2V$ gives a 6-symplectic
    structure on $V$, since the set of all degenerate two-forms on $V$ is a Pl\"ucker quadric, and all non-zero
    degenerate two-forms have two-dimensional kernel.

\hfill

There exist examples of non-degenerate $k$-symplectic structures for every $k$, as follows from
the construction described below.

Recall the definition of Clifford algebras (for the proofs of all statements about Clifford algebras, see \cite{ABS}).
We will always assume that all vector spaces are
either real or complex, however the definition below is valid for any field.
Consider a vector space $E$ with quadratic form $q\in S^2E^*$.
By definition, the Clifford algebra $\Cl(E,q)$ is the quotient of the tensor algebra $T^\bullet E$ by
the two-sided ideal generated by all tensors of the form $x\otimes x - q(x,x)\cdot 1$,
for all $x \in E$. The elements generating the ideal belong to the even part of the tensor algebra,
so the Clifford algebra is $\bbZ/2\bbZ$-graded, $\Cl(E,q) = \Cl^0(E,q)\oplus \Cl^1(E,q)$.
There exists an automorphism $\tau\colon \Cl(E,q)\to \Cl(E,q)$, which acts as the identity
on $\Cl^0(E,q)$ and as multiplication by $-1$ on $\Cl^1(E,q)$. There is also an
antiautomorphism of transposition which maps $x = x_1\cdot\ldots\cdot x_k$ to $x^t = x_k\cdot\ldots\cdot x_1$
for any $x_i\in E$. We will use the notation $\bar{x} = \tau(x^t)$ for $x\in\Cl(E,q)$.

One can check that the space $E$ is embedded into $\Cl(E,q)$ via the map $E=T^1E \to \Cl(E,q)$
induced by the projection $T^\bullet E\to \Cl(E,q)$.
The Clifford algebra has the following universal property. Let $\mathcal{A}$ be any associative
algebra with unit and $\alpha\colon E\to \mathcal{A}$ a map with the property that
$\alpha(x)^2 = q(x,x)\cdot 1_\mathcal{A}$ for all $x\in E$.
Then $\alpha$ can be uniquely extended to a morphism of algebras $\alpha'\colon \Cl(E,q)\to \mathcal{A}$
such that $\alpha'|_E = \alpha$.

Recall the definition of the group $\Pin(E,q)$ (see \cite{ABS}). We will denote by
$\Cl^\times(E,q)$ the multiplicative group of invertible elements.
Define 
$$\Pin(E,q) = \{x\in \Cl^\times(E,q) \mid \tau(x) E x^{-1} = E,\, x\bar{x} = 1\}.$$
One can check that $\Pin(E, q)$ acts on $E$ preserving the quadratic form $q$ and
is actually a double covering of the orthogonal group $\mathrm{O}(E,q)$.
The group $\Spin(E,q)$ is defined as the subgroup of even elements in $\Pin(E,q)$.

\hfill

\example
For each integer $k > 0$ we will construct an example of a $k$-symplectic structure on some vector space.
Start from a real $k$-dimensional vector space $E_\bbR$ with a negative-definite quadratic form $q$.
Consider the Clifford algebra $\Cl(E_\bbR,q)$ and the corresponding group $\Pin(E_\bbR,q)$.
This group is compact (because $q$ is negative-definite) and it acts on $E_\bbR$ by orthogonal transformations.

\hfill

\lemma\label{lem_embedding}
Consider a non-trivial real representation $\rho\colon \Cl(E_\bbR,q)\to \emrp(V_\bbR)$ in some real vector
space $V_\bbR$, $\dim V_\bbR = 4n$ (by this we mean a representation of $\Cl(E_\bbR,q)$ as an algebra
with unit). 
For any representation $V_\bbR$ of $\Cl(E_\bbR,q)$ as above, we have a $\Spin(E_\bbR, q)$-equivariant
embedding $$\alpha_\bbR\colon E_\bbR\hookrightarrow \Lambda^2 V^*_\bbR.$$
All non-zero two-forms in the image of $\alpha_\bbR$ are non-degenerate.

\hfill

\proof
If we have a representation $V_\bbR$, then the compact group $\Pin(E_\bbR,q)$ acts on $V_\bbR$,
and there exists a positive-definite invariant quadratic form $g\in (S^2V_\bbR^*)^{\Pin(E_\bbR,q)}$.

Note that any element $e\in E_\bbR$ with $e^2 = -1$ lies in $\Pin(E_\bbR,q)$. This follows from
direct computation: we have $\tau(e) = -e$, $e\bar{e} = 1$ and for any $x\in E_\bbR$ we have $\tau(e)x e^{-1} =
exe = 2 q(x,e)e - xe^2 = 2q(x, e)e + x \in E_\bbR$. From this and the invariance of $g$ with
respect to $\Pin(E_\bbR, q)$ we conclude that $g(e\cdot u,e\cdot v) = g(u,v)$ for any $u,v\in V_\bbR$.

We can define the two-form $\omega_e\in \Lambda^2V_\bbR^*$ by the formula $\omega_e(u,v) = g(e\cdot u,v)$
for $u,v\in V_\bbR$. This really defines an element in $\Lambda^2V_\bbR^*$ since
$g(e\cdot u,v) = g(e^2\cdot u, e\cdot v) = -g(u, e\cdot v) = -g(e\cdot v, u)$.

The quadratic form $q$ was chosen negative-definite, so every element $x\in E_\bbR$
is of the form $x = \xi e$ with $\xi\in \bbR$ and $e^2=-1$, so $\omega_x=\xi\omega_e\in \Lambda^2V_\bbR^*$.
We have constructed a map $\alpha_\bbR\colon E_\bbR\hookrightarrow \Lambda^2V_\bbR^*$, $x\mapsto \omega_x$,
which is clearly an embedding.

Note that this embedding is a morphism of $\Spin(E_\bbR,q)$-modules:
for any $h\in \Spin(E_\bbR,q)$ we have $\tau(h) = h$,
$\alpha_\bbR(h\cdot x) = \alpha_\bbR(hxh^{-1}) = \omega_{hxh^{-1}}$ and 
$\omega_{hxh^{-1}}(u,v) = g(hxh^{-1}\cdot u,v) = g(xh^{-1}\cdot u,h^{-1}\cdot v) = \omega_x(h^{-1}\cdot u, h^{-1}\cdot v)$,
so $\alpha_\bbR(h\cdot x) = h\cdot \omega_x$.

For all non-zero $x\in E_\bbR$ the form $\omega_x$ is non-degenerate, because the quadratic
form $g$ is chosen positive-definite (and in particular non-degenerate).
\endproof

\hfill

Consider the complexifications: $E=E_\bbR\otimes\bbC$, $\Cl(E,q)=\Cl(E_\bbR,q)\otimes\bbC$ and
$V = V_\bbR\otimes \bbC$ with the corresponding complex representation of $\Cl(E,q)$.
Then we get an embedding of complex vector spaces $\alpha\colon E\hookrightarrow \Lambda^2V^*$. 

So we have an irreducible $\Spin(E,q)$-module $E$ embedded into $\Lambda^2V^*$. Fix a volume form and
identify $\Lambda^{4n}V^*$ with $\bbC$. Then the wedge product
in $\Lambda^\bullet V^*$ induces a $\Spin(E,q)$-invariant polynomial $p\in S^{2n}E^*$,
$p(x)= \omega_x^{\wedge 2n}$ on $E$. This polynomial is non-zero, because the image
of $\alpha$ contains non-degenerate two-forms by 
\ref{lem_embedding}.

The group $\Spin(E,q)$ is a double covering of $\mathrm{SO}(E,q)$,
so $p$ is $\mathrm{SO}(E,q)$-invariant and by classical invariant theory this
polynomial has to be proportional to the $n$-th power of $q$.

We can define $\bOmega$ to be the image of $E$ in $\Lambda^2V^*$ under the map $\alpha$.
We have already seen that the set of degenerate
two-forms is a quadric. To see that $\bOmega$ is a $k$-symplectic structure it only
remains to check that all degenerate
two-forms in $\bOmega$ have kernels of dimension $2n$. This is proved in the following lemma
which we will also need later on.

\hfill

\lemma\label{dim_ker_lemma}
Let $V$ be a complex vector space of dimension $4n$. Let $\bOmega\subset\Lambda^2V^*$ be a
$k$-dimensional subspace. Identifying $\Lambda^{4n}_\C V^*$ with $\bbC$, consider a polynomial
$p\in S^{2n}{\bOmega}^*$ given by $p(\omega) = \omega^{\wedge 2n}$. Suppose that there exists
a non-degenerate quadratic form $q\in S^2{\bOmega}^*$, such that $p = q^n$.
Then all degenerate two-forms in $\bOmega$ have rank $2n$, in particular $\bOmega$
is a $k$-symplectic structure on $V$.

\hfill

\proof
Let $\omega_1\in \bOmega$ be such that $q(\omega_1,\omega_1)=0$, and let $\dim(\ker\omega_1) = 2r$.
Since $q$ is non-degenerate we can find $\omega_2\in\bOmega$ with $q(\omega_2,\omega_2)\neq 0$
and $q(\omega_1,\omega_2)\neq 0$. Consider the polynomial
$\tilde{p}(t) = p(\omega_1+t\omega_2)$. We have $\tilde{p}(t) = q(\omega_1+t\omega_2,\omega_1+t\omega_2)^n$,
so $\tilde{p}$ must have zero of order $n$ at $t=0$. But
$$
\tilde{p}(t) = (\omega_1+t\omega_2)^{\wedge2n} = t^r\omega_1^{\wedge(2n-r)}\wedge\omega_2^{\wedge r} + t^{r+1}\omega_1^{\wedge(2n-r-1)}\wedge\omega_2^{\wedge(r+1)}+\ldots+t^{2n}\omega_2^{\wedge 2n},
$$
and $\omega_1^{\wedge(2n-r)}\wedge\omega_2^{\wedge r}\neq 0$ because $\omega_2$ is non-degenerate.
So the order of zero at $t=0$ is $r$, hence $r=n$.
\endproof

\hfill

We can bound from below the dimension 
of a vector space carrying a $k$-symplectic structure.
In the proposition below we use the following terminology:

\hfill

\definition
A $k$-symplectic structure $(\bOmega, q)$ on a vector space $V$ is called
{\bf real} if $V$ is a complexification of a real vector space $V_\bbR$, $\bOmega$
is a complexification of a real subspace $\bOmega_\bbR$ in
$\Lambda^2V^*_\bbR$, and $q$ is a
real quadratic form.

\hfill

In particular, the $k$-symplectic structures from the previous example are real.
We will denote by $\Cl_{r,s}$ the Clifford algebra for a real $(r+s)$-dimensional vector
space with quadratic form of signature $(r,s)$, meaning that it has $r$ minuses and $s$ pluses.

\hfill

\proposition\label{_Clifford_from_k-symple_Proposition_}
\begin{enumerate}[label=(\alph*)]
\item
Let $(V, \bOmega, q)$ be a  vector space
equipped with a non-degenerate $q$-symplectic structure. 
Then $V$ is a $\Cl^0(\bOmega, q)$-module.
\item
Let $(V, \bOmega, q)$ be a vector space with a non-degenerate
real $k$-symplectic structure
where the real quadratic form $q$ has signature $(r,s)$. Then
the corresponding real vector space $V_\bbR$ 
has a structure of $\Cl_{r-1,s}$-module if $r>0$,
and of a $\Cl_{s-1, r}$-module if $s>0$.
\end{enumerate}

\proof
Consider a pair of elements $\omega_1, \omega_2\in \bOmega$, such that 
$q(\omega_1,\omega_1) = -1$ and $q(\omega_1,\omega_2) = 0$. These elements define linear maps
$\omega_i\colon V\to V^*$, and by definition of a $k$-symplectic
structure the map $\omega_1$ is an isomorphism.
So we can define the endomorphism $A = \omega_1^{-1}\omega_2\in \emrp(V)$. We claim that
\begin{equation}\label{eqn_a}
A^2 = q(\omega_2,\omega_2)Id.
\end{equation}
To prove this we need to find the eigenvalues of $A$: the operator
$A-\lambda Id = \omega_1^{-1}(\omega_2 - \lambda\omega_1)$
is degenerate if and only if the form $\omega_2 - \lambda\omega_1$ has non-trivial kernel.
By definition of a $k$-symplectic structure, this condition is equivalent to
$q(\omega_2-\lambda\omega_1, \omega_2-\lambda\omega_1) = 0$, hence the eigenvalues of $A$
are $\lambda_{\pm} = \pm q(\omega_2,\omega_2)^{1/2}$ and the claim follows.

Next fix an element $\omega_1\in \bOmega$ with $q(\omega_1,\omega_1) = -1$ and consider its $q$-orthogonal
complement $W = \{\omega_2\in \bOmega\mid q(\omega_1,\omega_2)=0\}$. Define a linear map
$\alpha\colon W\to \emrp(V)$, $\alpha(\omega_2) = \omega_1^{-1}\omega_2$. From the equation \eqref{eqn_a},
we see that $\alpha(\omega_2)^2 - q(\omega_2,\omega_2)Id = 0$. By the
universal property of the Clifford algebra, the map $\alpha$ can be extended to
$\alpha'\colon \Cl(W, q|_W)\to \emrp(V)$.

Observe that $\Cl(W, q|_W) \simeq \Cl^0(\bOmega,q)$.
This isomorphism can be constructed as follows. For any element $\eta\in \Cl(W, q|_W)$
consider the decomposition $\eta = \eta^0+\eta^1$ into even and odd parts. Then
the map $\eta\mapsto \eta^0 + \omega_1\eta^1\in \Cl^0(\bOmega,q)$ is the desired isomorphism,
which can be checked using the fact that $\omega_1^2 = -1$ and that $\omega_1$
commutes with even and anticommutes with odd elements from $\Cl(W,q|_W)$.
 This proves the first part of the proposition.

Next consider the case of a real $k$-symplectic structure with the quadratic form $q$ of
signature $(r,s)$. If $r>0$ we can choose a real element $\omega_1\in \bOmega_\bbR$
with $q(\omega_1,\omega_1) = -1$. Its real $q$-orthogonal complement $W_\bbR$ has a quadratic form
of signature $(r-1,s)$, and we obtain a representation of $\Cl_{r-1,s}$ in $V_\bbR$.

If $s>0$ we can choose a real element $\omega_1\in \bOmega_\bbR$
with $q(\omega_1,\omega_1) = 1$. In this case for any $\omega_2$ with $q(\omega_1,\omega_2)=0$
we have an operator $A = \omega_1^{-1}\omega_2$ with $A^2 = -q(\omega_2,\omega_2)Id$.
Then we obtain a representation of $\Cl(W_\bbR,q')$ in $V_\bbR$, where $W_\bbR$ is
the $q$-orthogonal complement to $\omega_1$ as above, and $q'=-q|_{W_\bbR}$. The form
$q'$ has signature $(s-1,r)$, so the second claim of the proposition follows.
\endproof

\hfill

\corollary\label{_Clifford_from_k-symple_Corollary_}
Let $(V, \bOmega, q)$ be vector space
equipped with a non-degegenerate $k$-symplectic structure. Then
$$
\dim_\bbC V = 2^{\left\lfloor(k-1)/2\right\rfloor}m
$$
for some positive integer $m$.

\hfill

\proof
By the previous proposition, $V$ is a non-trivial $\Cl^0(\bOmega,q)$-module, so we have a map
$\alpha'\colon \Cl^0(\bOmega,q)\to \emrp(V)$. The map $\alpha'$ is non-zero, so its image is isomorphic
to a quotient of $\Cl^0(\bOmega,q)$ by a proper two-sided ideal.
But $\Cl^0(\bOmega,q)$ is either the matrix algebra $\mathrm{Mat}({2^{(k-1)/2}},\bbC)$ if $k$ is odd
or the sum of two copies of the matrix algebra $\mathrm{Mat}({2^{k/2-1}},\bbC)$ if $k$ is even.
In any case, $V$ is a direct sum of
$m$ copies of the standard representation of the matrix algebra for some $m$. Hence the desired equality for $\dim V$.
\endproof

\hfill

\remark
Both the construction from the previous example and the proof of the previous proposition
involve Clifford algebras. Note, however, that these two constructions are not inverse to each
other. Not every $k$-symplectic structure arises from a real representation of a Clifford algebra
associated with a $k$-dimensional vector space. One example when this is not the case was already
mentioned above: this is a 6-symplectic structure on a 4-dimensional vector space $V$. In this case,
the proof of the previous proposition gives us a representation $\alpha'\colon \Cl(W,q)\to \emrp(V)$
with 5-dimensional vector space $W$ and $\Cl(W,q)$ isomorphic to $\mathrm{Mat}(4,\bbC)\oplus\mathrm{Mat}(4,\bbC)$.
The map $\alpha'$ here is not injective. Note that by dimension reasons $V$ cannot be a
representation of a Clifford algebra associated to a 6-dimensional vector space.
If we consider $V$ as a complexification of a real 4-dimensional vector space, then the corresponding
quadratic form on $\Lambda^2V^*_{\bbR}$ will be of signature $(3,3)$, and by the previous
proposition we have a $\Cl_{2,3}$-module structure on $V_\bbR$. The algebra $\Cl_{2,3}$ is isomorphic
to $\mathrm{Mat}(4,\bbR)\oplus\mathrm{Mat}(4,\bbR)$.


\subsection{Applications to hyperk\"ahler geometry}
\label{_k_symple_appli_Subsection_}

Next we will use our observations about $k$-symplectic structures to investigate 
submanifolds in a very general irreducible holomorphic
symplectic (IHS) manifold (see \ref{_IHS_Definition_}).

Consider an IHS manifold $X$ of dimension $2n$ with symplectic form $\sigma$.
We will call $X$ very general if it represents a point in the moduli space which lies in
the complement to a countable union of proper analytic subvarieties. One can prove that
a very general $X$ cannot contain submanifolds of odd dimension. Moreover, if there exists
a submanifold of dimension $2m$ inside such $X$, and $\gamma\in H^{4(n-m)}(X,\bbZ)$ is its fundamental
class, then $\gamma$ stays of Hodge type $(2n-2m,2n-2m)$
on any small deformation of $X$ (such a
subvariety is called {\bf absolutely trianalytic}, 
see \ref{_abs_triana_Definition_}).
The proof of these well-known facts can 
be found for example in \cite{GHJ}, section 26.3.

\hfill

\proposition\label{_abs_tria_torus_k_sy_Proposition_}
Let $X$ be an IHS manifold, $k=b_2(X)$ 
its second Betti number, and $T$ a compact complex torus
of dimension $2m$ immersed into $X$. Assume that $T$ is
absolutely trianalytic (this would follow if $X$ is
sufficiently general). Then $H_1(T,\bbC)$ carries
a non-degenerate real $k$-symplectic structure. 
The corresponding quadratic form has
signature $(k-3,3)$.

\hfill

\proof
Denote by $j\colon T\to X$ the immersion. Let $V = H_1(T,\bbC)$, then $\Lambda^2V^* = H^2(T,\bbC)$
and we have the restriction map $j^*\colon H^2(X,\bbC)\to H^2(T,\bbC)$. Fujiki relations (\ref{_pairing_subvari_Equation_})
imply that $j^*$ is injective: for any $\alpha, \beta\in H^2(X, \bbC)$ we have
$$
(j^*\alpha)^{2m-1}\cdot j^*\beta = \gamma\cdot \alpha^{2m-1}\cdot \beta = c_\gamma q(\alpha,\alpha)^{m-1}q(\alpha,\beta),
$$
and if $j^*\beta = 0$ then we must have $q(\alpha,\beta) = 0$ for any $\alpha$ with $q(\alpha,\alpha)\neq 0$.
Taking $\alpha$ to be a K\"ahler form we see that $\beta = 0$.

We define $\bOmega = \mathrm{Im}(j^*)\subset \Lambda^2V^*$. Then by Fujiki relations (\ref{_pairing_subvari_Equation_})
we have $(j^*\alpha)^{2m} = c_\gamma q(\alpha,\alpha)^m$ and by 
\ref{dim_ker_lemma} $\bOmega$
is a $k$-symplectic structure on $V$. This $k$-symplectic structure is real, because the
Beauville-Bogomolov form is real (the real structure on
$V$ is the standard real structure on homology).
The signature of Beauville-Bogomolov form is known to be $(b_2(X)-3, 3)$.
\endproof

\hfill

\corollary\label{_b_1_torus_Corollary_}
If a torus $T$ is immersed into a very general IHS manifold $X$ then
$$\dim_\bbC T\geqslant 2^{\left\lfloor(b_2(X)-1)/2\right\rfloor-1}.$$

\hfill

\proof
This follows directly from \ref{_abs_tria_torus_k_sy_Proposition_} and \ref{_Clifford_from_k-symple_Corollary_}
(note that the complex dimension of a torus is twice its first Betti number).
\endproof

\hfill

\remark\label{_10_dim_subva_Remark_}
In particular, this shows that a very 
general deformation of a 10-dimensional IHS manifold
$M$ constructed by O'Grady can not contain a 
complex torus, because it is known that in
this case $b_2(M)= 24$,
and this would imply $\dim_\bbC T\geqslant 2^{10}$, 
which is impossible.

\hfill

\remark\label{_6_dim_subva_Remark_}
For the 6-dimensional IHS manifold of O'Grady's, we have
$b_2 = 8$ and our estimate gives
$\dim_\bbC T\geqslant 4$. But in this case we can use the
fact that the real homology $H_1(T,\bbR)$ must carry a 
real 8-symplectic structure with quadratic
form of signature $(5,3)$. In this
case $H_1(T,\bbR)$ has to be a $\Cl_{4,3}$-module (and a $\Cl_{2,5}$-module).
The algebras $\Cl_{4,3}$ and $\Cl_{2,5}$ are isomorphic to $\mathrm{Mat}(8,\bbC)$
as real algebras. The minimal non-trivial real representation of $\mathrm{Mat}(8,\bbC)$ is
$\bbC^8\simeq \bbR^{16}$. This means that 
$H_1(T,\bbR)$ has to be at least 16-dimensional
which gives a contradiction.

\hfill

\section{Further remarks and open questions}

\subsection{$k$-symplectic structures on manifolds}

The notion of $k$-symplectic structure generalizes
that of a trisymplectic structure, which was originally defined in
\cite{_Jardim_Verbitsky:trisymple_}. 
Jardim and Verbitsky define a trisymplectic structure on a complex manifold $M$
to be a triple $\langle \Omega_1, \Omega_2, \Omega_3\rangle$
of symplectic forms satisfying the following condition.
Let $\bOmega$ be the vector space generated by
$\langle \Omega_1, \Omega_2, \Omega_3\rangle$, and
$Q\subset \bOmega$ the set of degenerate forms.
Then $Q$ is a non-degenerate quadric hypersurface,
all $v\in Q$ have constant rank, and all non-zero 
$v\in Q$ have rank $\frac 1 2 \dim M$.

This notion was studied in \cite{_Jardim_Verbitsky:trisymple_}
when $M$ is a complex manifold, and 
$\Omega_1, \Omega_2, \Omega_3$ holomorphic symplectic forms.
In this case $M$ admits an action of $\Mat(2, \C)$ in
its tangent space, preserving the space $\bOmega$,
and a torsion-free holomorphic connection preserving
the $\bOmega$ and $\Mat(2, \C)$-action. Trisymplectic
structures arise naturally in connection with the
mathematical instantons on $\C P^3$; indeed, the
moduli space of framed mathematical instantons on
$\C P^3$ is trisymplectic.

Trisymplectic structure is a special case of the
3-web, explored in \cite{_JV:Instantons_}, where it was studied
using the geometric approach coming back to Chern's doctoral
thesis \cite{_Chern:3-webs_}. A 3-web is a triple of 
involutive sub-bundles $S_1, S_2, S_3 \subset TM$
such that $TM= S_i \oplus S_j$ for any $i\neq j$.
Chern has constructed
a canonical connection on a manifold equipped with a
3-web; for 3-webs arising from trisymplectic structures,
this connection turns out to be torsion-free.

Trisymplectic manifold can be characterized in terms of
holonomy of this connection. These are manifolds equipped
with a torsion free connection $\nabla$ such that its
holonomy group $\Hol(\nabla)$ is contained in a group 
$G=SL(2, \C)\cdot Sp(2n, \C)$ acting on $4n$-dimensional
complex space $\C^{4n}=\C^2\otimes \C^{2n}$, with
$SL(2, \C)$ acting tautologically on the first tensor
factor, and $Sp(2n, \C)$ on the second tensor factor.

Trisymplectic structures occur naturally in hyperk\"ahler
geometry, for the following reason. Let $M$ be a
hyperk\"ahler manifold, and $\Tw(M)$ its twistor space,
that is, a total space of its twistor family 
(\ref{_twistor_family_Definition_}).
Twistor space is a complex manifold, but it is 
non-algebraic and non-K\"ahler when $M$ is compact.

In \cite{_HKLR_} (see also \cite{_Simpson:hyperka-defi_}  
and \cite{_Verbitsky:desing_}) the twistor 
space was used to define (possibly singular)
hyperk\"ahler varieties. It turns out that
a component $\Sec(M)$ of the moduli of rational
curves on $\Tw(M)$ is equipped with a real structure,
in such a way that a connected component $\Sec_\R(M)$ of its
set of real points is identified with $M$.
Then one can describe the hyperk\"ahler structure
on $M$ as a certain geometric structure on $\Sec_\R(M)$.

In \cite{_Jardim_Verbitsky:trisymple_}
it was shown that $\Sec(M)$ is equipped with a 
trisymplectic structure, and its restriction to
$M=\Sec_\R(M)$ gives the triple of symplectic forms
$\omega_I, \omega_J, \omega_K$. One can understand this
result by identifying $\Sec(M)$ and a complexification of
$M$. Then the trisymplectic structure on $\Sec(M)$ 
is obtained as a complexification of the triple
$\omega_I, \omega_J, \omega_K$. 




A related notion of even
Clifford structures on Riemannian manifolds was 
introduced by A.~Moroianu
and U.~Semmelmann in \cite{_MS_}. An {\bf even Clifford structure} on
a Riemannian manifold $M$ is a sub-bundle $A\subset \End(TM)$
which is closed under
multiplication and fiberwise isomorphic to an even Clifford
algebra acting on the tangent bundle $TM$ orthogonally.
An even Clifford structure is called {\bf parallel}
if $A$ is preserved by the Levi-Civita connection.
If $A$ is trivial, this Clifford structure is called
{\bf flat}. Moroianu and Semmelmann classified all
manifolds admitting a parallel Clifford structure.

If we are
given a $k$-symplectic structure, then it induces an
action of an even Clifford algebra, 
as one can see from \ref{_Clifford_from_k-symple_Proposition_}.
However, the metric induced by symplectic forms and
the Clifford algebra action is not necessarily positive definite. 
If one removes the positive definiteness condition, then 
any $k$-symplectic structure induces an even, flat Clifford structure,
in the sense of Moroianu-Semmelmann. 

For $k=3$, the corresponding flat Clifford structure is 
also preserved by the Levi-Civita connection, hence parallel
(\cite{_Jardim_Verbitsky:trisymple_}). It is not clear whether
the same is true for $k>3$ (see \ref{_parallel_dim_3_Question_} below).

In \cite{_Moroianu_Pilca_}, A. Moroianu and M. Pilca 
obtain topological restrictions on manifolds admitting
Clifford structures. They prove that a manifold with 
non-vanishing Euler characteristic cannot admit a Clifford structure
with rank bigger than 16.

\subsection{Open questions and possible directions of research}

The notion of a $k$-symplectic manifold is certainly
an intriguing one, and still very little understood.

Let $\bOmega$ be a $k$-symplectic structure on $M$.
A general $l$-dimensional subspace $\bOmega'\subset \bOmega$
obviously gives an $l$-symplectic structure on $M$. 
This means, in particular, 
that any $k$-symplectic manifold is equipped
with a $3(k-3)$-dimensional family $R$ of
3-symplectic structures. As shown in 
\cite{_Jardim_Verbitsky:trisymple_}, a 3-symplectic
manifold is equipped with a canonical torsion-free
connection preserving the 3-symplectic structure
(this connection is called {\bf the Chern connection};
not to be confused with the Chern connection defined
on holomorphic Hermitian bundles).
Given a 3-symplectic structure $r\in R$, denote
the corresponding Chern connection by $\nabla_r$.

\hfill

\question
Are all $\nabla_r$ equal for all $r\in R$? If so, what is
the holonomy of the Chern connection associated with a
$k$-symplectic structure this way?

\hfill

\question\label{_parallel_dim_3_Question_}
Let $M$ be a complex manifold equipped with a real
structure $\iota$, and $\bOmega$ an $\iota$-invariant 
complex $k$-symplectic structure. The real part
$\bOmega_\R$ restricted to the real part $M^\iota$
gives a $k$-dimensional space of symplectic forms of
constant rank. Degenerate forms $\omega\in \bOmega$ 
correspond to the real points in the quadric $Q$.
However, if the set of real points of $Q$ is 
empty, all forms in 
$\bOmega_\R \subset \Lambda^2_\R(M^\iota)$
are non-degenerate.
When $k=3$, $\bOmega$ defines a hyperk\"ahler
structure on $M^\iota$, and any hyperk\"ahler
structure can be defined this way. What happens
if we start with $k$-symplectic structure?
Presumably, we obtain a generalization of
hyperk\"ahler structures, with the quaternion
algebra replaced by a bigger-dimensional
Clifford algebra. Such a manifold 
would admit a very special family of
hyperk\"ahler structures, parametrized
by a Grassmanian of 3-dimensional subspaces
in a $k$-dimensional space.

\hfill

In \cite{_Jardim_Verbitsky:trisymple_}, 
a geometric reduction construction was established
for trisymplectic structures. 
Trisymplectic quotient was defined in such a way that
a trisymplectic quotient of a trisymplectic manifold
$M$ equipped with an action of $d$-dimensional
compact Lie group $G$ preserving the
trisymplectic structure is a $(\dim M-4d)$-dimensional
complex manifold equipped with a trisymplectic structure.
This construction can be considered as a complexification
of the hyperk\"ahler reduction of \cite{_HKLR_}.

\hfill

\question
Is there a geometric reduction construction
for $k$-symplectic structures?

\hfill

3-symplectic geometrical reductions are 
known in complex (\cite{_Jardim_Verbitsky:trisymple_})
and real versions (\cite{_Dancer_Swann:hypersy_quotie_}).

{\bf Acknowledgement:}
We are grateful to K. Oguiso for interesting and fruitful
discussions and to F. Bogomolov, D. Kaledin and A. Kuznetsov for useful
comments.

{\small
\noindent \sc
Andrey Soldatenkov\\
Institut f\"ur Mathematik\\
Humboldt-Universit\"at zu Berlin\\
Unter den Linden 6\\
10099 Berlin\\
{\tt aosoldatenkov@gmail.com}\\
\ \  \\
\noindent {\sc Misha Verbitsky\\
{\sc Instituto Nacional de Matem\'atica Pura e
              Aplicada (IMPA) \\ Estrada Dona Castorina, 110\\
Jardim Bot\^anico, CEP 22460-320\\
Rio de Janeiro, RJ - Brasil }\\
also:\\
{\sc Laboratory of Algebraic Geometry,\\
National Research University HSE,\\
Department of Mathematics, 7 Vavilova Str. Moscow, Russia,}\\
\tt  verbit@impa.ru}.
}

\end{document}